\documentclass[11pt]{article}
\usepackage[top=1in, bottom=1in, left=1in, right=1in]{geometry}
\usepackage{paper}
\usepackage{nicefrac}

\usepackage{todonotes} 
\usepackage{comment}

\usepackage[title]{appendix}

\newcommand{\Ntrain}{N_\text{train}}
\newcommand{\Ntest}{N_\text{test}}

\newcommand{\stf}{solve-training framework}
\newcommand{\st}{solve-training}
\newcommand{\Stf}{Solve-training framework}

\newcommand{\ftf}{fit-training framework}
\newcommand{\ft}{fit-training}
\newcommand{\Ftf}{Fit-training framework}

\renewcommand{\vec}[1]{\boldsymbol{#1}}

\begin{document}

\title{Variational training of neural network approximations of
solution maps for physical models}

\author{
    Yingzhou Li$^{\,\dagger}$,
    Jianfeng Lu$^{\,\dagger\,\ddagger}$,
    Anqi Mao$^{\,\sharp}$
  \vspace{0.1in}\\
  $\dagger$ Department of Mathematics, Duke University\\
  $\ddagger$ Department of Chemistry and Department of Physics,
  Duke University\\
  $\sharp$ Department of Mathematics, Shanghai Jiao Tong University\\
}

\maketitle

\begin{abstract} 
A novel solve-training framework is proposed to train neural network
in representing low dimensional solution maps of physical models.
Solve-training framework uses the neural network as the ansatz of the
solution map and trains the network variationally via loss functions
from the underlying physical models. Solve-training framework avoids
expensive data preparation in the traditional supervised training
procedure, which prepares labels for input data, and still achieves
effective representation of the solution map adapted to the input
data distribution. The efficiency of solve-training framework
is demonstrated through obtaining solution maps for linear and
nonlinear elliptic equations, and maps from potentials to ground
states of linear and nonlinear Schr\"odinger equations.
\end{abstract}

{\bf Keywords.} Neural network; solution map for PDEs; fast algorithm;
hierarchical matrix; unsupervised training.

\section{Introduction}
\label{sec:intro}

Simulation of physical models has been one of main driven forces for
scientific computing.  Physical phenomena at different scales, e.g.,
macroscopic scale, microscopic scale, etc., are characterized by
Newton's laws of motion, Darcy's law, Maxwell's equations,
Schr\"odinger equation, etc. Solving these equations efficiently,
especially those nonlinear ones, has challenged computational
scientists for decades and led to remarkable development in algorithms
and in computing hardware. As the rise of machine learning,
particularly deep learning, many researchers have been attempting to
adopt artificial neural networks (NN) to represent the
high-dimensional solutions or the low-dimensional solution maps. This
paper proposes a variational training framework for solving the
solution map of low-dimensional physical models via NNs. Here we
emphasize solving a solution map in contrast with fitting a solution
map, where solving can be to some extent viewed as unsupervised
learning with input functions only and fitting refers to supervised
learning with both input functions and the corresponding solutions.

Solving the solution map for physical models is feasible due to an
intrinsic difference between the physical problems and other
data-driven problems, e.g., handwriting recognition, speech
recognition, spam detection, etc. Indeed, for physical models, the
solution maps are governed by well-received equations, which are often
expressed in partial differential equations (PDEs), whereas the
conventional machine learning tasks such as image classification rely
on human labeled data set without explicit expression for the
underlying model.  Benefiting from such a difference, we design loss
functions based on the PDEs, in another word, we adopt the model
information into the loss functions, and solve the solution map
directly without knowing solution functions.

\subsection{Related work}

A number of recent work utilized NNs to address physical models.
Generally, they can be organized into three groups: representing
solutions via NNs, representing solution maps via NNs, and optimizing
traditional iterative solvers via NNs.  Representing solutions of
physical models, especially high-dimensional ones, has been a
long-standing computational challenge.  NN with multiple input and
single output can be used as an ansatz for the solutions of physical
models or PDEs, which is first explored in \cite{Lagaris1998} for
low-dimensional solutions. Many high-dimensional problems, e.g.,
interacting spin models, high-dimensional committor functions, etc.,
have been recently considered for solutions using NN ansatz with
variants optimization strategies~\cite{Berg2018, Carleo2017, E2018,
  Han2018a, Khoo2019, Sirignano2018, Rudd2015}. NN, in this case, is
valuable in its flexibility and richness in representing
high-dimensional functions.

Representing the solution map of a nonlinear problem is challenging as
well.  For linear problems, the solution map can be represented by a
simple matrix (\textit{i.e.}, Green's function for PDE
problems). While the efficient representation for solution map is
unknown for most nonlinear problems. Traditional methods in turn solve
nonlinear problem via iterative methods, e.g., fixed point
iteration. Since NN is able to represent high dimensional nonlinear
mappings, it has also been explored in recent literature to represent
solution maps of low-dimensional problems on mesh grid, see
e.g,~\cite{Fan2019a, Fan2018, Fan2019, Han2017, Khoo2017, Khoo2018,
  Li2018, Long2018, Sun2003, Tang2017}. These NNs are fitted by a set
of training data with solution ready, i.e., labeled data.  Most works
from the first two groups focus on creative design of NN
architectures, in particular trying to incorporate knowledge of the
PDE into the representation.

The last group, very different from previous two, adopts NN to
optimize traditional iterative methods~\cite{Greenfeld2019, Hsieh2019,
  Katrutsa2017, Mishra2018}. Once the iterative methods are optimized
on a set of problems, generalization to different boundary conditions,
domain geometries, and other similar models, is explored and can be
sometimes guaranteed~\cite{Hsieh2019}.

\subsection{Main Idea}
\label{sec:idea}

The goal of this work is to propose a new paradigm of training neural
networks to approximate solution maps for physical PDE models, which
does not rely on existing PDE solvers or collected solution data.  The
main advantages come in two folds: the new training framework removes
the expensive data preparation cost and obtains an input-data-adaptive
NN with better accuracy in terms of intrinsic criteria from the PDE
after training.

\begin{figure}[htp]\centering
    \begin{subfigure}{0.48\textwidth}
        \includegraphics[width=\hsize]{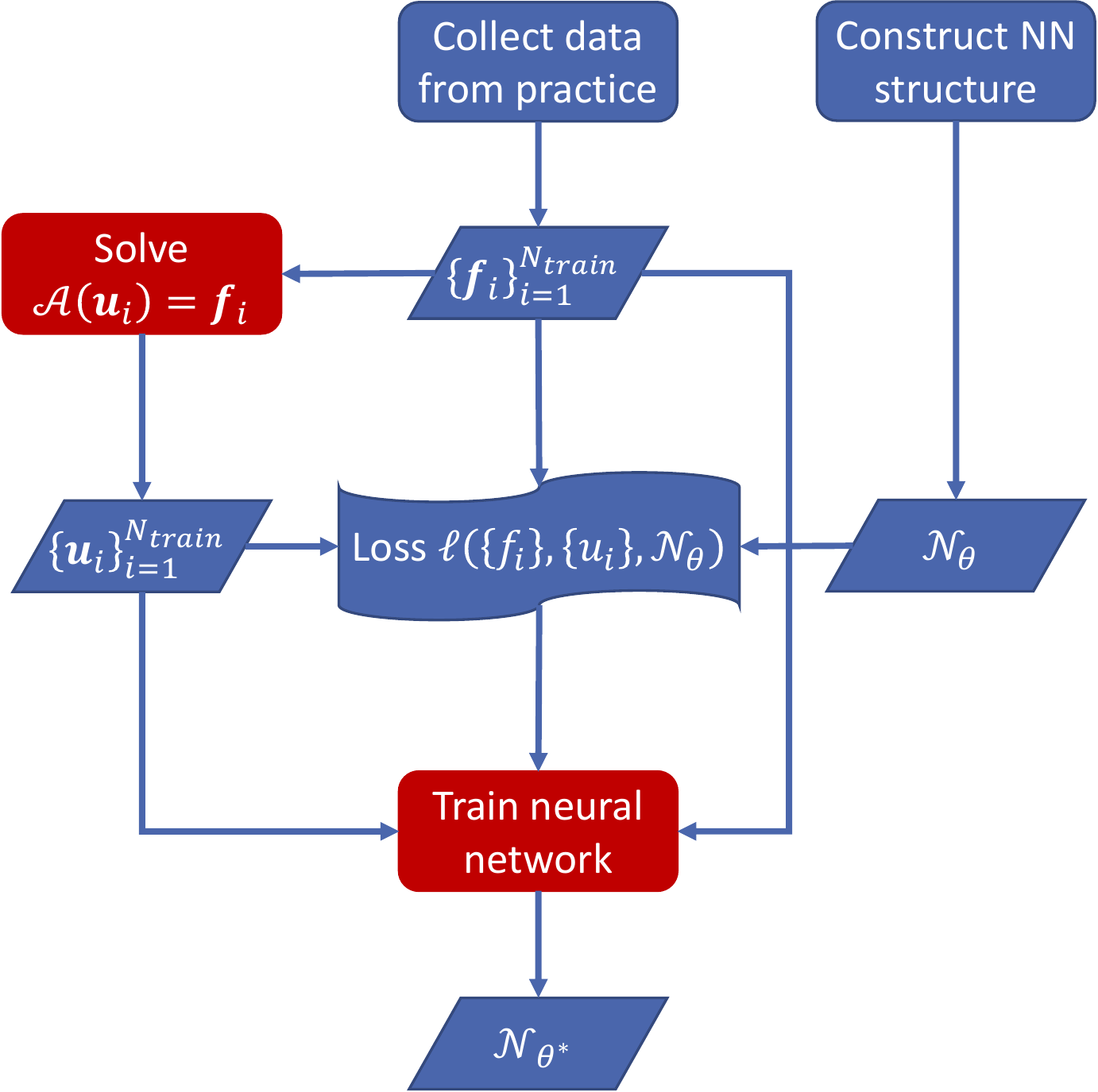}
        \caption{\Ftf{}.}
    \end{subfigure}    
    \hspace{0.4cm}
    \begin{subfigure}{0.48\textwidth}
        \includegraphics[width=\hsize]{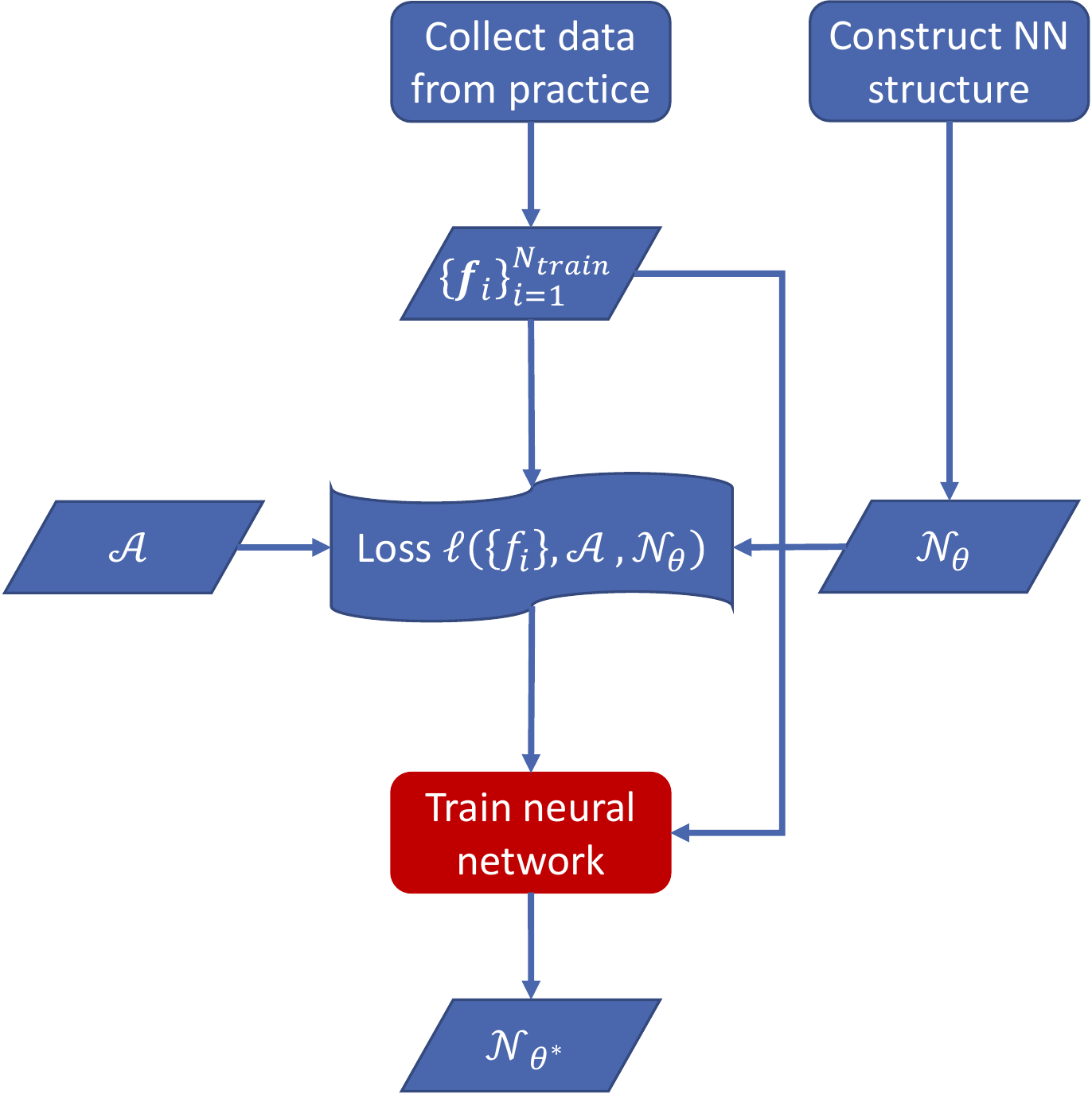}
        \caption{\Stf{}.}
    \end{subfigure}    
    \caption{Flowchart of (a) \ftf{} and (b) \stf{}. Rounded corner
    rectangles in general indicate processes whereas the red ones are
    computationally expensive processes; parallelograms indicate data
    or NNs; punched tapes are loss functions which require creative
    design depending on the problem; line arrows indicate dependencies
    between blocks.} \label{fig:flowchart}
\end{figure}

We now explain the main idea of the new training framework through the
example of solving a (possibly nonlinear) system of equations. Later in
Section~\ref{sec:solveeig}, we will show that the training framework
can be applied to solve the solution map of linear and nonlinear
eigenvalue problem as well.

Let us consider a system of equations, written as
\begin{equation}
    \calA(\vec{u}) = \vec{f},
\end{equation}
with $\vec{u}$ and $\vec{f}$ denoting solution and input functions
on a mesh and $\calA$ denoting a discretized forward operator. The
goal here is to obtain the solution map, i.e.,
\begin{equation}
    \vec{u} = \calA^{-1}(\vec{f}) \approx \calN_\theta(\vec{f}),
\end{equation}
which is approximated by a NN $\calN_\theta$ parameterized by $\theta$.
The input data $\vec{f}$ is usually collected from practice following
an unknown distribution $\calD_f$, denoted as $\vec{f} \sim \calD_f$.
Ideally, NN not only approximates the inverse map $\calA^{-1}$,
but also adapts to the distribution, i.e., $\calN_\theta \approx
\calA^{-1}\mid_{\calD_f}$.

Almost all previous works design NN $\calN_\theta$ based on properties
of the problem and then fit the solution map following the flowchart
in Figure~\ref{fig:flowchart}~(a). We call such a training procedure
\emph{\ftf{}}. In practice, there are two procedures for generating
training data pair, $\{\vec{u}_i,\vec{f}_i\}_{i=1}^{\Ntrain}$
\begin{enumerate}
    \item[(TD.1)] \label{TD1} Collect a set of
    $\{f_i\}_{i=1}^{\Ntrain}$ and solve $\vec{u}_i =
    \calA^{-1}(\vec{f}_i)$ with traditional methods;
    \item[(TD.2)] \label{TD2} Randomly generate a set of solution
    data, $\{\vec{u}_i\}_{i=1}^{\Ntrain}$ and evaluate $\vec{f}_i =
    \calA(\vec{u}_i)$.
\end{enumerate}
The first procedure suffers from expensive traditional method in
solving $\calA^{-1}$, especially when the equation is complicated and
nonlinear. While, the resulting training data set from the first
procedure follows the practical distribution $\calD_f$. Hence fitting
$\calN_\theta$ with this data set approximates
$\calA^{-1} \mid_{\calD_f}$.  The second procedure is efficient in
generating data since $\calA$ is usually cheap to evaluate. However,
the data set lacks proper distribution, i.e.,
$\vec{f}_i \not\sim \calD_f$. An accurate solution map requires
$\calN_\theta$ to be a good representation of $\calA^{-1}$ instead of
$\calA^{-1} \mid_{\calD_f}$, which is much more difficult to
approximate in general.

Under our new training framework, illustrated in
Figure~\ref{fig:flowchart}~(b), $\vec{u}_i$ is not required in the
loss function hence not required in the training procedure. Instead,
the forward mapping $\calA$ is brought into the loss function. One
simplest example of such a loss function in the sense of mean square
error is
\begin{equation} \label{eq:loss}
    \ell \Big( \{\vec{f}_i\}_{i=1}^{\Ntrain}, \calA, \calN_\theta \Big)
    = \frac{1}{\Ntrain} \sum_{i=1}^{\Ntrain} \norm{\vec{f}_i - \calA
    \big( \calN_\theta (\vec{f}_i) \big)}^2,
\end{equation}
where $\calA$ in our implementation is represented by an NN with fixed
parameters. In contrast to the \ftf{}, the proposed training framework
-- \emph{\stf{}} -- has at least three advantages:
\begin{enumerate}
    \item Solving $\vec{u}_i = \calA^{-1}(\vec{f}_i)$ via expensive
        traditional method is not needed;
    \item The trained NN $\calN_\theta$ is able to capture $\calA^{-1}
        \mid_{\calD_f}$;
    \item The parameters obtained through \stf{}
        minimizes the ``$\calA$-norm'' between $\calN_\theta(\vec{f})$
        and $\vec{u}_{*} = \calA^{-1}(\vec{f})$, i.e.,
        \begin{equation}
            \norm{\vec{f} - \calA \big( \calN_\theta(\vec{f} ) \big)} =
            \norm{\vec{u}_* - \calN_\theta(f)}_{\calA},
        \end{equation}
        if $\calA$ satisfies assumptions such that
        $\norm{\cdot}_{\calA}$ is well-defined.
\end{enumerate}
Regarding the last point above, \ftf{} minimizes 2-norm between
$\calN_\theta(\vec{f})$ and $\vec{u}_*$, which corresponds to least
square fitting for linear operators.  Hence we claim that \stf{} is
more likely to obtain an NN $\calN_\theta$ which solves $\calA(\vec{u})
= \vec{f}$ given $\vec{f} \sim \calD_f$. Other than the neural network
approximation error, \stf{} contains one more source of approximation,
the discretization error of the real forward operator to $\calA$. Since
the discretized forward map $\calA$ is represented by a fixed neural
network, which does not increase the number of trainable parameters in
the training part, we can use high-accuracy discretization schemes to
significantly reduce the discretization error and make the error much
less than the neural network approximation error. The discretization
error is also contained in the \ftf{}, since we need to solve the
discretized equation to provide training data.

In this work, we demonstrate the power of the \stf{} through training the
NNs representing the solution maps of linear and nonlinear systems and
linear and nonlinear eigenvalue problems. We remark that while finishing
the work, we discovered some very recent works~\cite{Bar2019, Zhu2019}
aiming at solving inverse problems, whose training strategy shares some
similarity with the solve-training framework we proposed above.

\subsection{Organization}

The rest of the paper is organized as the following.
Section~\ref{sec:solvesystem} applies the \stf{} to solving linear and
nonlinear systems. The corresponding numerical results are attached
right after the problem description. Similar structure applies to
Section~\ref{sec:solveeig}, in which we solve linear and nonlinear
eigenvalue problems rising from Schr\"odinger equations.  Finally,
Section~\ref{sec:conclusion} concludes the paper with discussions on
extensibility.

\section{Solving linear and nonlinear systems}
\label{sec:solvesystem}

This section aims to show that the \stf{} can be applied to obtain
the NN representation of the solution maps of linear and nonlinear
systems. The main idea of \stf{} for solving systems has been illustrated
in Section~\ref{sec:idea}. We will demonstrate the efficiency of the
\stf{} through two examples, linear elliptic equation and nonlinear
elliptic equation.

\subsection{Linear elliptic equations}
\label{sec:solvelinearsystem}

In this section, we focus on the two dimensional linear variable
coefficient elliptic equations with periodic boundary condition, i.e.,
\begin{equation} \label{eq:linearellipticeq}
    -\grad \cdot a(x,y) \grad u(x,y) = f(x,y), \quad (x,y) \in \Omega =
    [0,1)^2,
\end{equation}
where $a(x,y)>0$ denote variable coefficients. Such an equation
appears in a wide range of physical models governed by Laplace's
equation, Stokes equation, etc. For \eqref{eq:linearellipticeq}
with constant coefficients, the inverse operator~\footnote{When
\eqref{eq:linearellipticeq} has constant coefficient with periodic
boundary condition, the most efficient method should be fast Fourier
transform.} has an explicit Green's function representation and can
be applied efficiently with quasilinear cost through fast multipole
methods~\cite{Fong2009, Greengard1987, Ying2004}, or other related
fast algorithms~\cite{Corona2015, Ho2016a}.  When the coefficient
is variable, then the operator in \eqref{eq:linearellipticeq} is
discretized into a sparse matrix, and solved via iterative methods
with efficient preconditioners~\cite{Briggs2000, Hackbusch1999,
Hackbusch2002, Ho2016, Li2017a, Xia2009}. Among these preconditioners,
$\calH$-matrices~\cite{Hackbusch1999, Hackbusch2002} are efficient
preconditioners of simplest algebraic form and the structures with
modifications are recently extended to NN structures~\cite{Fan2019a,
Fan2018}. While, the construction of the $\calH$-matrices for
the inverse of variable coefficient elliptic equations requires
sophisticated matrix-vector multiplication on structured random
vectors~\cite{Lin2011}. In this section, we adopt the original structure
of $\calH$-matrix in $\calN_\theta$ with and without ReLU layers. \Stf{}
then provides a method to construct the $\calH$-matrix with a limited
number of input functions.

The discretization of \eqref{eq:linearellipticeq} used here is the
five-point stencil on a $64 \times 64$ uniform grid. The discretization
points are $\{x_i,y_j\}_{i,j=0}^{63}$ with $x_i = i/64$ and $y_j =
j/64$. And the discrete variable coefficient $a(x_i,y_j)$ is a Chess
board field as,
\begin{equation} \label{eq:a}
    a(x_i,y_i) =
    \begin{cases}
        10, & \floor{ \frac{i+j}{8} } \equiv 0\, (\mathrm{mod} \, 2) \\
        1, & \floor{ \frac{i+j}{8} } \equiv 1\, (\mathrm{mod} \, 2) \\
    \end{cases}.
\end{equation}
And we generate $\Ntrain$ random vectors $\{\vec{f}_i\}_{i=1}^{\Ntrain}$
as the training data. Each $\vec{f}_i$ is a vector of length $64^2$ with
each entry being uniform random on $\Big[ -\sqrt{3/64^2}, \sqrt{3/64^2}
\Big]$ such that $\bbE \Big(\norm{\vec{f}_i}\Big) = 1$~\footnote{Notice
that normalization here is not important for the linear model and
we will use relative error as the measure in the later numerical
results. However, the NN package Tensorflow~\cite{tensorflow2015} uses
float32 as the default data format and such a normalization reduces
the impact of numerical errors.} and subtract its mean to incorporate
with the periodic boundary condition. This procedure defines $\calD_f$,
which will be less emphasized for linear model in this section. Another
set of $\Ntest = 5,000$ random vectors of the same distribution,
$\{\vec{g}_i\}_{i=1}^{\Ntest}$, is generated for testing purpose. The
reported relative error is calculated as follows,
\begin{equation} \label{eq:relerr}
    \frac{1}{\Ntest} \sum_{i=1}^{\Ntest} \frac{\norm{\vec{g}_i - \calA
    \Big( \calN_\theta(\vec{g}_i) \Big)}}{\norm{\vec{g}_i}}.
\end{equation}

Four $\calH$-matrices are generated and compared in this
section. The structures of all these $\calH$-matrices are
generated from bi-partition of the domain up to four layers and
each low-rank submatrix is of rank 96. Readers are referred to
the textbook~\cite{Hackbusch2015} for the detailed structure
of an $\calH$-matrix. The first $\calH$-matrix is constructed
directly from the inversion of the discretized sparse matrix and
each low-rank block is constructed via the truncated singular value
decomposition (SVD).  This $\calH$-matrix is close-to-optimal in the
standard $\calH$-matrix literature and is used as the baseline for
the comparison. We denote it as \emph{$\calH$-matrix~(SVD)} in the
later content. The second and third $\calH$-matrices are constructed
in the same way in Tensorflow~\cite{tensorflow2015}. The second one
is initialized with random coefficients and then trained, whereas
the third one is initialized with the baseline $\calH$-matrix and
then trained. They are denoted as \emph{NN-$\calH$-matrix~(rand
init)} and \emph{NN-$\calH$-matrix~(SVD init)} respectively. The
last $\calH$-matrix uses the same structure but with each small
dense block coupled with $5$ ReLU layers in the similar fashion
as in~\cite{Fan2018}. This $\calH$-matrix is initialized with SVD
coefficients and the ReLU part is initialized in a way such that the
initial output (no train) is the same as that of $\calH$-matrix~(SVD)
and then trained. It is denoted as \emph{NLNN-$\calH$-matrix~(SVD
init)}.

We train the later three $\calH$-matrices under the \stf{}
with Adam optimizer~\cite{Kingma2015}. The batch size is
100 for all trainings. For NN-$\calH$-matrix~(SVD init) and
NLNN-$\calH$-matrix~(SVD init), a fixed stepsize $2 \times 10^{-6}$
is used. While, for NN-$\calH$-matrix~(rand init), the stepsize
is initialized as $2 \times 10^{-4}$ following a steady exponential
decay to $2 \times 10^{-6}$. For each $\calH$-matrix, we train the NN
for three times and report the best among them.  Default values are
used for all other unspecified hyperparameters.

\subsubsection*{Numerical Results}

We first compare the performance of the first three $\calH$-matrices
described above through numerical experiments under the \stf{}.

\begin{table}[htp]
    \centering
    \begin{tabular}{lccccc}
        \toprule
        $\calH$-matrix & \# Epoch & Train loss & Test loss
        & Train rel err & Test rel err \\
        \toprule
        $\calH$-matrix~(SVD)
        & 0    & 2.21e-3 &2.20e-3  & 4.68e-2 & 4.66e-2\\
        NN-$\calH$-matrix~(rand init)
        &25000 & 2.09e-4 & 3.38e-4 & 1.44e-2 & 1.83e-2\\
        NN-$\calH$-matrix~(SVD init)
        & 2000 & 2.43e-4 & 3.40e-4 & 1.55e-2 & 1.83e-2\\
        \bottomrule
    \end{tabular}
    \caption{Train and test relative error of \stf{}
    in linear elliptic equation for different $\calH$-matrices. The train 
    and test data sets are of size $\Ntrain=10000$ and $\Ntest=5000$.
    } \label{tab:lpde}
\end{table}

\begin{table}[htp]
    \centering
    \begin{tabular}{lccccc}
        \toprule
        $\calH$-matrix & \# Epoch & Train loss & Test loss
        & Train rel err & Test rel err \\
        \toprule
        $\calH$-matrix~(SVD)
        & 0    & 2.20e-3 &2.20e-3  & 4.66e-2 & 4.66e-2\\
        NN-$\calH$-matrix~(rand init)
        &60000 & 1.34e-4 & 4.18e-2 & 1.15e-2 & 1.96e-1\\
        NN-$\calH$-matrix~(SVD init)
        & 6000 & 1.58e-4 & 7.20e-4 & 1.25e-2 & 2.67e-2\\
        \bottomrule
    \end{tabular}
    \caption{Train and test relative error of \stf{} in linear
    elliptic equation for different $\calH$-matrices. The train and
    test data sets are of size $\Ntrain=4000$ and $\Ntest=5000$.}
    \label{tab:lpde_train4000}
\end{table}

\begin{figure}[htp]\centering
	\begin{subfigure}{0.495\textwidth}
		\includegraphics[width=\hsize]{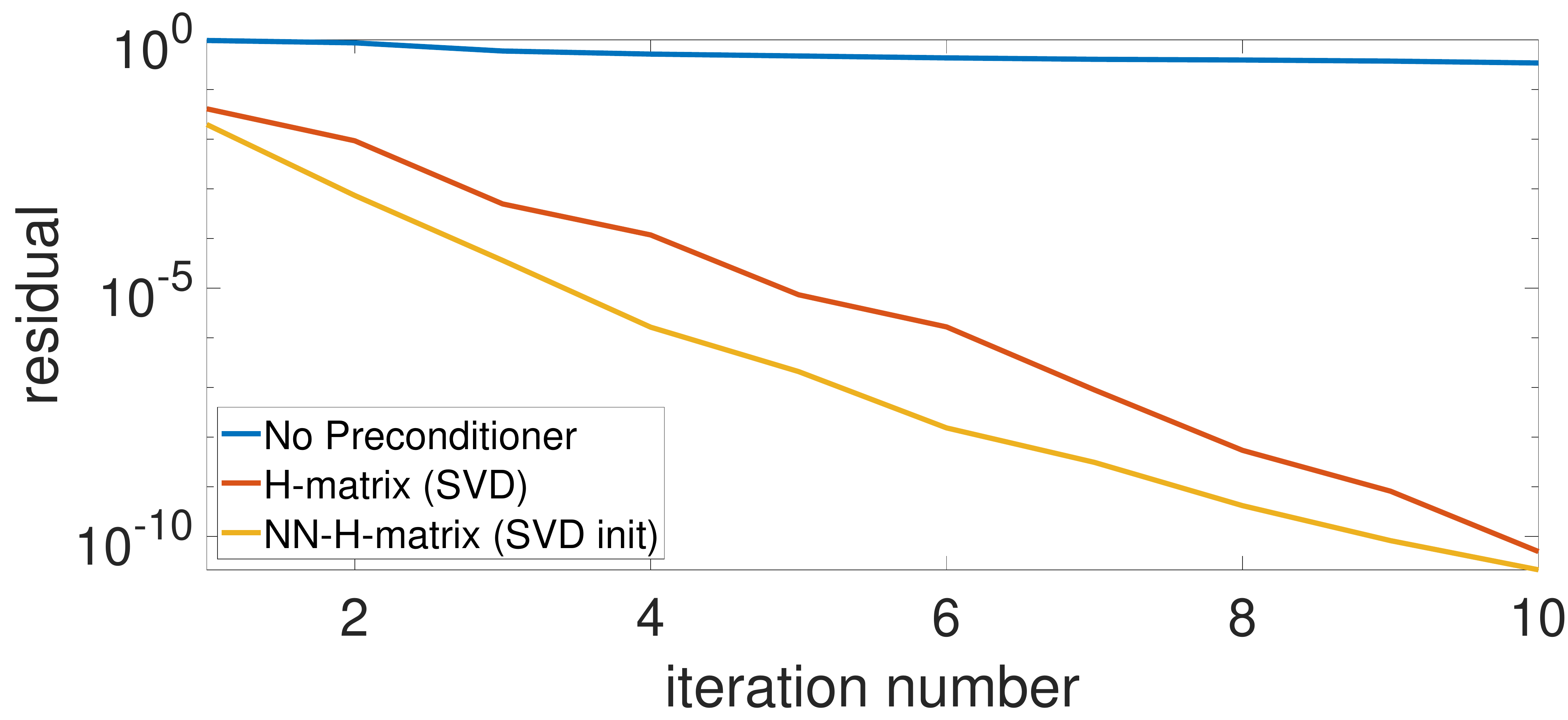}
		\caption{}
	\end{subfigure}    
	\begin{subfigure}{0.495\textwidth}
		\includegraphics[width=\hsize]{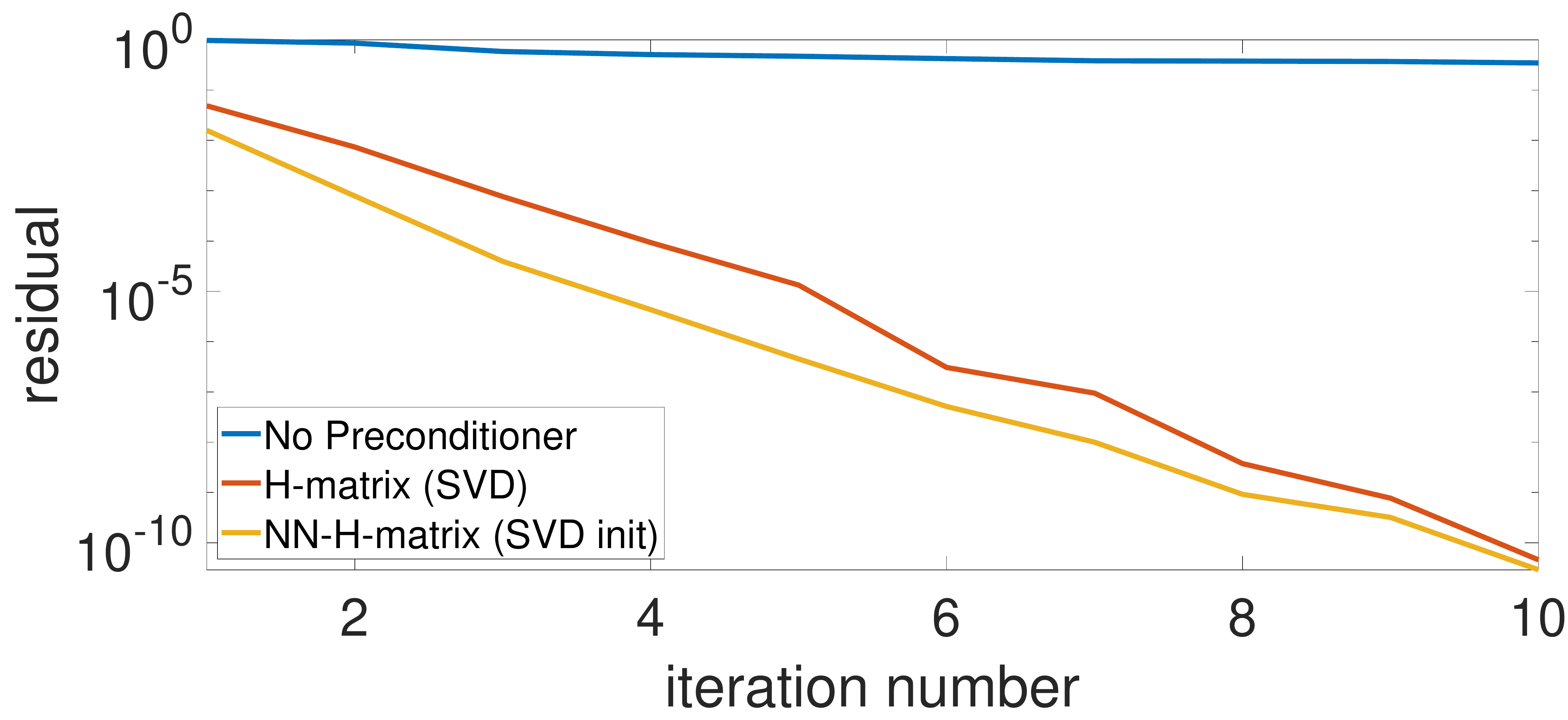}
		\caption{}
	\end{subfigure}    
	\caption{Examples of residual of using the conjugate gradient method to solve the linear elliptic equation, where $\calH$-matrix~(SVD) and the trained NN-$\calH$-matrix~(SVD init) are applied as preconditioners respectively. The train and test data sets are of size $\Ntrain=10000$ and $\Ntest=5000$.}
	\label{fig:precond}
\end{figure}

\begin{figure}[htp]
    \centering
    \includegraphics[width=0.6\textwidth]{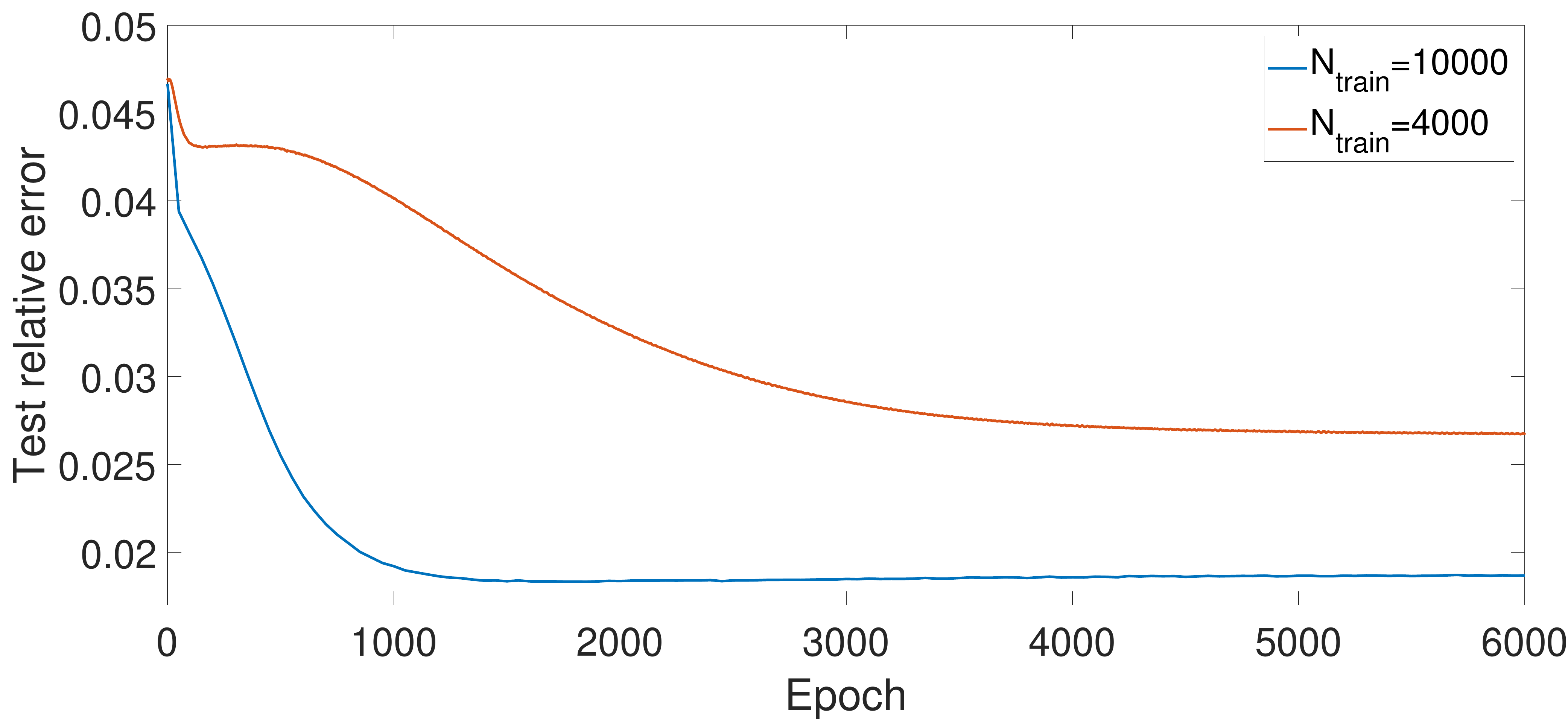}
    \caption{The test relative error of NN-$\calH$-matrix~(SVD
    init) against epochs.}
    \label{fig:H-mat3_validationerr}
\end{figure}

Table~\ref{tab:lpde} and Table~\ref{tab:lpde_train4000} present the
number of epoches, train loss, test loss, train relative error, and test
relative error for the first three $\calH$-matrices with $\Ntrain =
10000$ and $\Ntrain = 4000$ respectively. All of these matrices share
exactly the same structure and are all linear operators. Since these
$\calH$-matrices are trained on a uniform random input function and
they are linear, the test relative error is generalizable to other
non-normalized general input functions. When $\Ntrain = 10000$, we notice
that NN-$\calH$-matrix~(rand init) and NN-$\calH$-matrix~(SVD init)
achieve almost identical losses and relative errors after training under
the \stf{}, although the efficient training of NN-$\calH$-matrix~(rand
init) requires more aggressive choice of stepsize in the beginning of
the training.  Since we inject part of the information of the system
into the NN through the carefully designed architecture, training
under \stf{} is able to approximate the solution map with the number
of training data $\Ntrain$ smaller than the size of the matrix,
i.e., $\Ntrain = 4000$. In this case, NN-$\calH$-matrix~(SVD init) is
able to achieve similar results as that with $\Ntrain = 10000$. While
NN-$\calH$-matrix~(rand init) achieves similar train results but less
accurate test results.

In general, after training, the relative error for NN-$\calH$-matrices
is better than that of the $\calH$-matrix~(SVD), which means
that the low-rank approximation in $\calH$-matrix can be further
improved. Low-rank approximation through truncated SVD achieves
best 2/F-norm approximation locally in each block, whereas
the trained NN-$\calH$-matrix achieves near-optimal low-rank
approximation in the global sense. Also, Figure \ref{fig:precond}
shows examples of residual of using the conjugate gradient method
to solve the linear elliptic equation, where $\calH$-matrix~(SVD)
and the trained NN-$\calH$-matrix~(SVD init) are applied as
preconditioners respectively. The trained NN-$\calH$-matrix
achieves smaller residual than $\calH$-matrix after the same number
of iteration steps.  Hence the \stf{} can be applied to, either
obtaining the $\calH$-matrix representation of the inverse variable
coefficient elliptic operator, or further refine some existing fast
algorithms and achieves better approximation accuracy. In addition,
Figure~\ref{fig:H-mat3_validationerr} shows the refinement step is
quite efficient. The initial test relative error equals $0.0466$
and monotonically drops as the training goes on. After roughly $1000$
epoches, training with $\Ntrain = 10000$ samples, the test relative
error reaches a plateau with values about $0.02$.

\begin{figure}[htp]
    \centering
    \begin{subfigure}{0.495\textwidth}
        \includegraphics[width=\hsize]{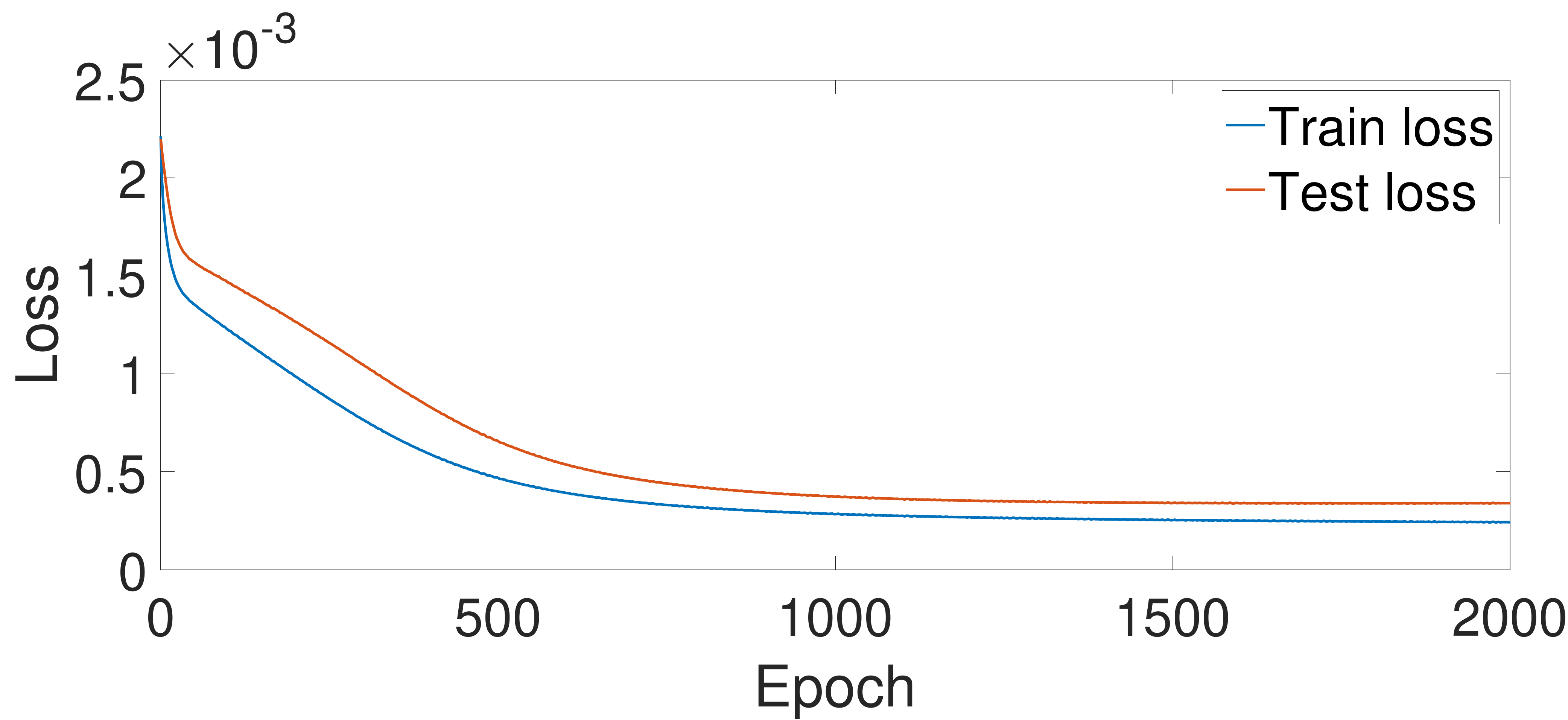}
        \caption{}
    \end{subfigure}    
    \begin{subfigure}{0.495\textwidth}
        \includegraphics[width=\hsize]{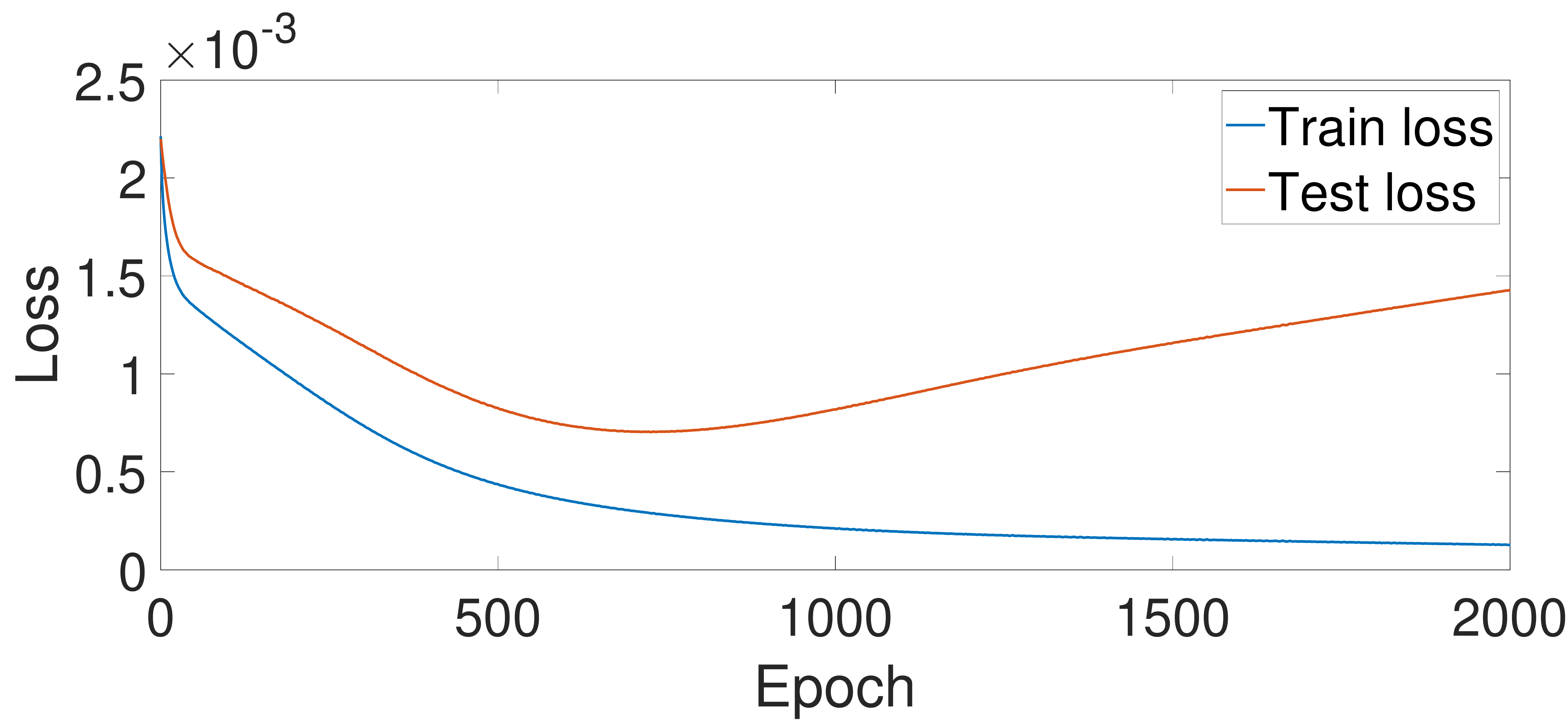}
        \caption{}
    \end{subfigure}    
    \caption{Train loss and test loss for (a) NN-$\calH$-matrix~(SVD
    init) and (b) NLNN-$\calH$-matrix~(SVD init) against epochs.}
    \label{fig:H-mat4_loss}
\end{figure}

NLNN-$\calH$-matrix~(SVD init) is a nonlinear operator approximating
the discrete inverse matrix of \eqref{eq:linearellipticeq}.
Figure~\ref{fig:H-mat4_loss}~(b) shows that its training behavior
and as a comparison, Figure~\ref{fig:H-mat4_loss}~(a) shows the
training behavior of NN-$\calH$-matrix~(SVD init). We observe severe
over fitting issue occurs in training NLNN-$\calH$-matrix~(SVD
init). Hence, for linear elliptic equations, training linear
operators under the \stf{} is a more preferred strategy to represent
$\calA^{-1}\mid_{\calD_f}$ accurately.

\subsection{Nonlinear elliptic equation}
\label{sec:solvenonlinearsystem}

In this section, we focus on a two dimensional nonlinear variable
coefficient elliptic equation with periodic boundary condition,
i.e.,
\begin{equation} \label{eq:nonlinearellipticeq}
        -\grad \cdot a(x,y) \grad u(x,y) + b u^3(x,y) =
        f(x,y), \quad (x,y) \in \Omega = [0,1)^2,
\end{equation}
where $a(x,y)>0$ denote variable coefficients and $b$ denotes the
strength of the nonlinearity. This equation adds a cubic nonlinear term
to \eqref{eq:linearellipticeq} with coefficient $b$ and is also
related to the nonlinear Sch\"ordinger equation introduced later in
Section~\ref{sec:nlschrodinger}. Solving such an equation can be
achieved through a carefully designed fixed-point iteration.  Hence
obtaining training data $\{\vec{u}_i, \vec{f}_i\}_{i=1}^{\Ntrain}$
for \ftf{} is expensive. In this section, we apply the \stf{} to train
an NN, $\calN_\theta$, to represent the solution map from $\vec{f}$
to $\vec{u}$.

We discretize \eqref{eq:nonlinearellipticeq} over a $32 \times 32$
uniform grid with five-point stencil. The discretization points are
$\{x_i,y_j\}_{i,j=0}^{31}$ with $x_i = i/32$ and $y_j = j/32$ and
the discrete variable coefficient $a(x_i,y_j)$ is defined the same
as in \eqref{eq:a}. And $b$ is set to be 0.1.

In order to show the advantage of the \stf{} regarding data
distributions, we train and test this example with three different
ways of generating training data $\{\vec{f}_i\}_{i=1}^{\Ntrain}$:
\begin{enumerate}
    \item[(D1)] Generate a set of solution data
    $\{\vec{u}_i\}_{i=1}^{\Ntrain}$ with each entry following
    the normal distribution $\calN(0,10^{-4})$, and then evaluate
    $\vec{f}_i = \calA(\vec{u}_i)$.
    \item[(D2)] Generate a set of solution data
    $\{\vec{u}_i\}_{i=1}^{\Ntrain}$ and each $\vec{u}_i$ is
    a convolution of a Gaussian kernel of standard deviation
    $\frac{1}{16}$ with a random vector with each entry following the
    normal distribution $\calN(0,10^{-4})$. Then evaluate $\vec{f}_i
    = \calA(\vec{u}_i)$.
    \item[(D3)] Generate $\{\vec{f}_i\}_{i=1}^{\Ntrain}$ and each
    $\vec{f}_i$ is a convolution of a Gaussian kernel of standard
    deviation $\frac{1}{16}$ with a random vector three entries of
    which are randomly picked up to follow the uniform distribution
    $\mathcal{U}(0.1,0.3)$ and other entries equal to 0. All
    $\vec{f}_i$ is subtracted by its mean and hence is mean zero.
\end{enumerate}
We assume D3 generates the input data, which is regarded as the
collected data. D1 and D2 are two designed distributions generating
training data for the purpose in \hyperref[TD2]{TD.2} and hence
the expensive traditional solving step is avoided.  For each kind
of data, $\Ntrain = 50,000$ training samples and $\Ntest = 5,000$
testing samples of the same distribution are generated. The reported
relative error is calculated as \eqref{eq:relerr}.

To authors' best knowledge, no existing NN structure is designed to
represent the solution map of \eqref{eq:nonlinearellipticeq}. Since the
focus of this paper is not on the creative design of the NN structure,
we construct a simple NN but by no means an efficient one for the
task. Thanks to the universal approximation theory~\cite{Barron1993},
the solution map can be represented by a single layer NN with accuracy
depending on the width. We construct an NN with one fully connected
layer of 10240 units using ReLU activation function to approximate the
solution map of \eqref{eq:nonlinearellipticeq}.

We train the single layer NN under the \stf{} and also the \ftf{}
for comparison purposes with Adam optimizer. The batch size is 100
for all trainings. The NN is trained for 10000 epochs with stepsize
$2 \times 10^{-4}$.  Default values are used for all other unspecified
hyperparameters.

\subsubsection*{Numerical results}

Since D3 is assumed to be the given data following the
distribution of interest $\calD_f$, to train an NN representing
$\calA^{-1}\mid_{\calD_f}$ under \ftf{}, one has to solve
$\calA^{-1}$ by expensive traditional methods.  Another choice
for \ftf{} is to obtain training data from other distributions and
generalize to $\calD_f$.  Hence we proposed D1 and D2 as alternative
choices of the distribution and validate the generalizability to
D3. However, \stf{} approximates $\calA^{-1}\mid_{\calD_f}$ directly.
If $\calA^{-1}\mid_{\calD_f}$ is much easier for NN to represent than
$\calA^{-1}$, then the generalizability of the trained NN under \stf{}
to D1 and D2 should be limited.
  
\begin{table}[htp]
    \centering
    \begin{tabular}{ccllcll}
        \toprule
        & Train Data & Train loss & Train rel err & Test Data
        & Test loss  & Test rel err \\
        \toprule
        \multirow{2}{*}{\ft{}}
        & D1 & 3.13e-6 & 6.22e-4 & D3 & -       & 1.60e+1 \\
        & D2 & 2.52e-5 & 1.65e-2 & D3 & -       & 1.29e+0 \\
        \midrule
        \multirow{3}{*}{\st{}} & 
        \multirow{3}{*}{D3} &
        \multirow{3}{*}{3.00e-4} &
        \multirow{3}{*}{1.71e-3}
        &       D1 & 1.81e+3 & 4.34e+0 \\
        & & & & D2 & 3.05e+1 & 9.90e-1 \\
        & & & & D3 & 6.18e-4 & 1.96e-3 \\
        \midrule
        \st{} & D3-N & 1.03e-2 & 1.00e-2 & D3 & 6.31e-4 & 2.01e-3 \\
        \bottomrule
    \end{tabular}
    \caption{Relative error of \stf{} and \ftf{} for the nonlinear
    elliptic equation given different kinds of train and
    test data.} \label{tab:nlpde}
\end{table}

Table~\ref{tab:nlpde} illustrates test relative error
of \ftf{} and \stf{} for the nonlinear elliptic equation
\eqref{eq:nonlinearellipticeq} given different choices of train and
test data.  Comparing the second to last row in Table~\ref{tab:nlpde}
against two rows of \ftf{}, we conclude that \stf{} successfully
trained a NN for approximating $\calA^{-1}\mid_{\calD_f}$ since both
the train and test relative error achieves almost three digits of
accuracy. While NN trained under \ftf{} on a synthetic distribution D1
and D2 achieves excellent relative error on training data but fails to
produce reliable prediction for data in D3. Since D2 is smoother than
D1, which has closer distribution to D3, NN trained under \ftf{} on
D2 performs sightly better than that on D1.  Here we also include the
test loss and relative error of NN, which trained under \stf{} on D3,
on D1 and D2 in Table~\ref{tab:nlpde}. The success of the approximation
of the solution map is distribution dependent. \Stf{} is also robust
to the noise in data. D3-N in the last row of Table \ref{tab:nlpde}
is obtained by adding Gaussian white noise $\mathcal{N}(0,0.01)$
to D3. Though the training loss and relative error increase, the
test loss and relative error on D3 are only slightly higher than the
results of NN trained directly on D3 without any noise.  Regarding the
computational cost of \ftf{} and \stf{}, although we have extra cost
in applying $\calA$ in the train procedure, it is negligible comparing
to the cost of other parts in NN. In practice, we observe that the
runtime for the train procedures of all experiments for both \ftf{}
and \stf{} in this section are about the same.

\section{Solving linear and nonlinear eigenvalue problem}
\label{sec:solveeig}

This section aims to show that the \stf{} not only can be applied to
solve linear and nonlinear systems but also can be applied to solve
the solution map of smallest eigenvalue problems.

Given an abstract eigenvalue problem as
\begin{equation} \label{eq:eig}
    \calA(u(x),V(x)) = E u(x),
\end{equation}
where $\calA$ denotes the operator, $V(x)$ denotes the external
potential, $u(x)$ is the eigenfunction corresponding to the eigenvalue
$E$. Many physical problems interest in the computation of the ground
state energy and ground state wavefunction, i.e., the smallest
eigenvalue and the corresponding eigenfunction. Since $V(x)$ is
the input external potential function, we define the solution
map of \eqref{eq:eig} being the map from $V(x)$ to $u(x)$
which corresponds to the smallest eigenpair, i.e., $\calM(V(x))
= u(x)$. In the discrete setting, we abuse notation $\calM$ to
represent the discrete solution map, i.e., $\calM(\vec{V}) = \vec{u}$,
where $\vec{V}$ and $\vec{u}$ denote discrete potential function and
ground state wavefunction respectively. Earlier work~\cite{Fan2018}
shows that a specially designed NN is able to capture the solution
map $\calM$ given the distribution of $\vec{V}$ as $\calD_V$, i.e.,
\begin{equation}
    \calN_\theta(\vec{V}) \approx \calM\mid_{\calD_V} (\vec{V})
    = \vec{u},
\end{equation}
where $\calN_\theta$ denotes an NN parameterized by $\theta$.

Under \ftf{}, as in most of previous work, the training of $\calN_\theta$
relies on the following loss function,
\begin{equation}\label{equ:snetloss_eig}
    \ell \Big(\{\vec{V}_i\}_{i=1}^{\Ntrain},
    \{\vec{u}_i\}_{i=1}^{\Ntrain}, \calN_\theta \Big) = \frac{1}{\Ntrain}
    \sum_{i=1}^{\Ntrain} \norm{\calN_\theta (\vec{V}_i) - \vec{u}_i}^2,
\end{equation}
where $\{\vec{u}_i,\vec{V}_i\}_{i=1}^{\Ntrain}$ are training data. For
eigenvalue problem, it becomes infeasible to obtain the training data
through the forward mapping of randomly generated $\vec{u}_i$, since
the ground state energy $E_i$ is unknown. Hence obtaining training data
requires solving a sequence of expensive linear/nonlinear eigenvalue
problems.

For the eigenvalue problems, it is computationally beneficial if the
training can be done under \stf{}. However, designing the loss function
is tricky and problem dependent under the \stf{}. Even if we assume $E$
is represented by another NN, the na\"ive loss function, i.e., two norm
of the difference of two sides of \eqref{eq:eig}, does not work, since
such a loss has multiple global minima corresponding to all eigenpairs
of \eqref{eq:eig}. Hence, solving the na\"ive loss function in many
cases does not give the eigenpair associated with smallest eigenvalue.
For the following linear Schr\"odinger and nonlinear Schr\"odinger
equations, we propose two loss functions to train the NN under \stf{}.

\subsection{Linear Schr\"odinger equation}
\label{sec:lschrodinger}

This section focuses on training the solution map of the smallest
eigenvalue problems of the linear one-dimensional Schr\"odinger
equation as,
\begin{equation}\label{eq:linearschrodingereq}
    \begin{split}
        & -\Delta u(x) + V(x) u(x) = E u(x), \quad x
        \in \Omega = [0,1)\\
        & \mathrm{s.t.} \int_{\Omega} u(x)^2 \diff x = 1, \quad
        \text{and} \int_{\Omega} u(x) \diff x > 0,
    \end{split} 
\end{equation}
with periodic boundary condition. The second positivity constraint
in \eqref{eq:linearschrodingereq} can be dropped since if $u(x)$
is the eigenvector associated to the smallest eigenvalue then so
is $-u(x)$. Besides, the right-hand side of the first constraint in
\eqref{eq:linearschrodingereq} can take any positive constant since
this eigenvalue problem is linear.

The external potential is randomly generated to simulate crystal
with two different atoms in each unit cell, i.e., $V(x)$ is randomly
generated via,
\begin{equation} \label{eq:potentialfun}
    V(x) = - \sum_{i=1}^{2} \sum_{j=-\infty}^{\infty}
    \frac{\rho^{(i)}}{\sqrt{2\pi T}} \exp(
    -\frac{\abs{x-j-c^{(i)}}^2}{2T}),
\end{equation}
where $c^{(i)}\sim \calU(0,1)$ for $i=1,2$ are the locations of two
atoms, and $\rho^{(i)} \sim \calU(10,40)$ and $T\sim \calU(2,4)\times
10^{-3}$ characterize the mass and electron charges of atoms.
Here $\calU(a, b)$ denotes the uniform distribution on the interval
$(a, b)$.

In this section, the linear Schr\"odinger equation
\eqref{eq:linearschrodingereq} is discretized on a uniform grid
in $[0,1)$ with $2048$ grid points.  The Laplace operator in
\eqref{eq:linearschrodingereq} is then discretized by the second-order
central difference scheme. Each input vector $\vec{V}$ composes of
the external potential $V(x)$ evaluated at grid points.

We propose a loss function as in the quadratic form,
\begin{equation} \label{eq:lschrodingerloss}
    \ell \Big(\{\vec{V}_i\}_{i=1}^{\Ntrain}, \calA,
    \calN_\theta \Big) = \sum_{i=1}^{\Ntrain} \ev**{-\Delta +
    \vec{V}_i}{\calN_\theta(\vec{V}_i)},
\end{equation}
which depends only on $\{\vec{V}_i\}_{i=1}^{\Ntrain}$, $\calA$,
$\calN_\theta$. When $\calN_\theta$ outputs a normalized result, each
term in the loss function is a variational form of the eigenvalue.
Hence, minimizing the loss function gives the ground state energy
$E$ if $\calN_\theta$ is able to capture the solution map of
\eqref{eq:linearschrodingereq}.

In addition to $\Ntrain$ training set $\{\vec{V}_i\}_{i=1}^{\Ntrain}$,
another set of $\Ntest$ random external potential
vectors of the same distribution as \eqref{eq:potentialfun},
$\{\vec{W}_i\}_{i=1}^{\Ntest}$, is generated for testing purpose. The
train and test loss as \eqref{eq:lschrodingerloss} is the summation
of all smallest eigenvalues and does not show the approximation
power of $\calN_\theta$ to the solution map given the distribution
of $\vec{V}$. Hence, we compare the output of trained $\calN_\theta$
against the underlying true smallest eigenvector and report the
relative error, which is calculated as follows,
\begin{equation} \label{eq:lschrodingerrelerr}
    \frac{1}{\Ntest} \sum_{i=1}^{\Ntest} \norm{ \vec{u}_i -
    \calN_\theta(\vec{W}_i)},
\end{equation}
where $\vec{u}_i$ is the normalized smallest eigenvector
corresponding to $\vec{W}_i$ for $i=1, \dots, \Ntest$. Equation
\eqref{eq:lschrodingerrelerr} is called the relative error
since $\vec{u}_i$ for $i=1,\dots,\Ntest$ are normalized, i.e.,
$\norm{\vec{u}_i} = 1$.

Since \citet{Fan2018} designed an $\calH$-matrix inspired NN structure,
called $\calH$-net in this paper, and successfully fitted the solution
map of nonlinear Schr\"odinger equations under the \ftf{}, we adopt their
structure here with a small modification to enforce the normalization
constraint. More precisely, the $\calH$-net is generated with eight
layers and each low-rank block is of rank $6$.
We vary the number of ReLu layers in the dense block, and the number
is denoted as $K$ in the later content. One extra normalization layer
is added in the end of $\calH$-net, i.e.,
\begin{equation}
    \calN_\theta(\vec{V}) =
    \frac{\widetilde{\calN}_\theta}{\bigl\lVert\widetilde{\calN}_\theta \bigr\rVert},
\end{equation}
where $\widetilde{\calN}_\theta$ is the regular
$\calH$-net~\cite{Fan2018} and $\calN_\theta$ is the NN used in this
section. Since the normalization layer does not involve any parameter,
the same $\theta$ is used for both $\widetilde{\calN}_\theta$ and
$\calN_\theta$.

We train $\calN_\theta$ under the \stf{} and also the \ftf{}
for comparison with Adam optimizer. The batch size is 100 for all
trainings. $\calN_\theta$ is trained for 60,000 epochs with stepsize
as $2 \times 10^{-4}$. Default values are used for all other unspecified
hyperparameters.

\subsubsection*{Numerical Results}

We first compare the performance of $\calN_\theta$ trained under \stf{}
for different number of train data set size $\Ntrain$ and different
number of ReLU layers $K$ through numerical experiments.

\begin{table}[htp]
    \centering
    \begin{tabular}{cccc}
        \toprule
        $\Ntrain$ & $\Ntest$ & Train rel err & Test rel err \\
        \toprule
        500       & 5000     & 9.46e-2       & 1.01e-1 \\
        1000      & 5000     & 2.07e-2       & 2.54e-2 \\
        5000      & 5000     & 7.11e-3       & 8.16e-3 \\
        20000     & 20000    & 7.84e-3       & 8.15e-3 \\
        \bottomrule
    \end{tabular}
    \caption{Relative error of $\calN_\theta$ with $K=5$ trained
    under \stf{} for linear Schr\"odinger equation given different
    sizes of train and test data set.} \label{tab:lse_size}
\end{table}

Table~\ref{tab:lse_size} presents the relative errors for different
$\Ntrain$ and $\Ntest$ with $K=5$. The test relative error decreases
as $\Ntrain$ increases. However, $\Ntrain = 5000$ train samples have
already been able to provide near-optimal results, since both the
train relative error and the test relative error stay similar for
$\Ntrain = 5000$ train samples and $\Ntrain = 20000$ train samples.
Hence, in this section, we adopt $\Ntrain = 5000$ and $\Ntest = 5000$
for all later experiments.

\begin{table}[htp]
    \centering
    \begin{tabular}{cccc}
        \toprule
        $K$ & $N_{params}$ & Train rel err & Test rel err \\
        \toprule
        1   & 15184        & 1.70e-1       & 1.72e-1  \\
        3   & 34236        & 2.88e-2       & 3.00e-2  \\
        5   & 57156        & 7.11e-3       & 8.16e-3  \\
        7   & 83944        & 5.87e-3       & 7.59e-3  \\
        \bottomrule
    \end{tabular}
    \caption{Relative error of $\calN_\theta$ trained under \stf{}
    for linear Schr\"odinger equation with different number of ReLU
    layers $K$. The train and test data sets are of size $\Ntrain=5000$
    and $\Ntest=5000$.} \label{tab:lse_K}
\end{table}

\begin{figure}[htp]\centering
    \begin{subfigure}{0.495\textwidth}
        \includegraphics[width=\hsize]{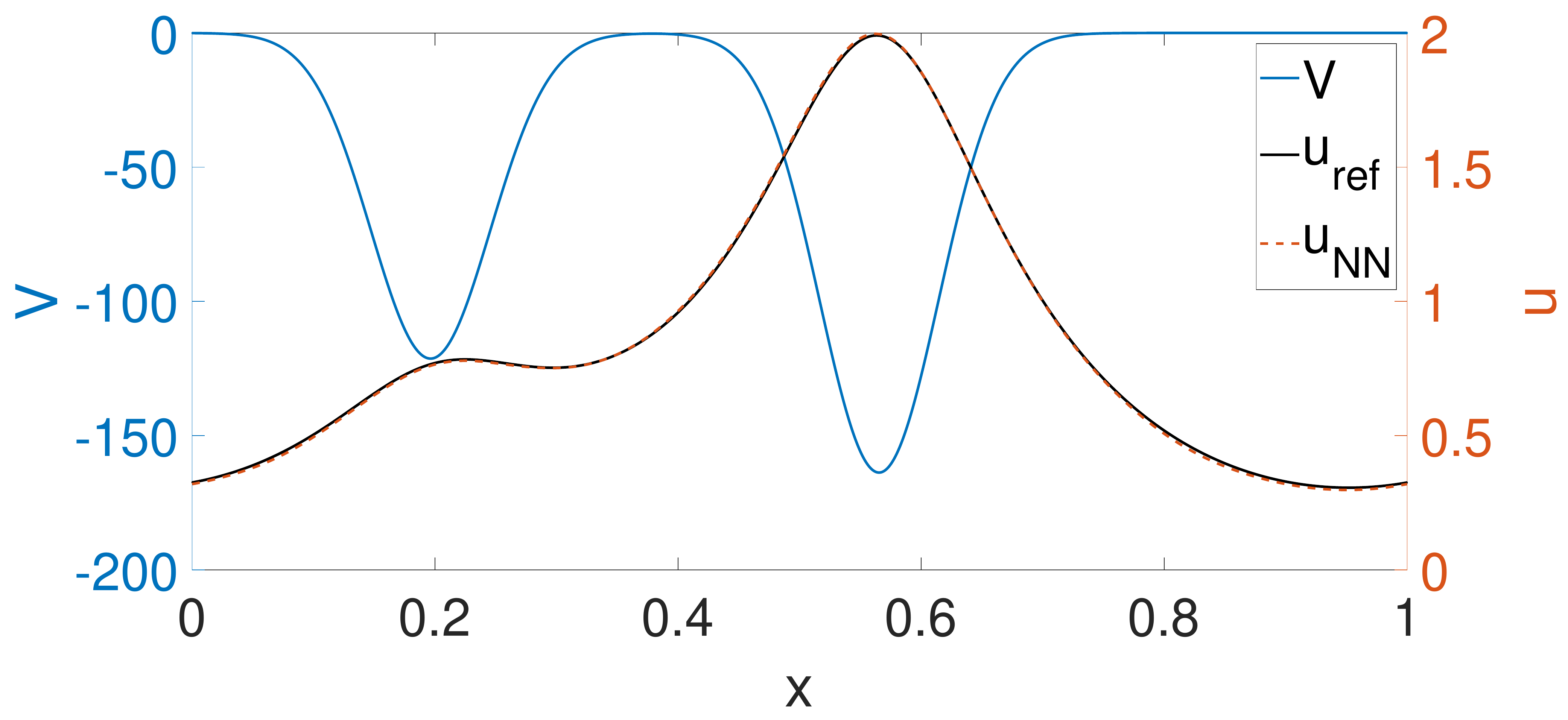}
        \caption{}
    \end{subfigure}    
    \begin{subfigure}{0.495\textwidth}
        \includegraphics[width=\hsize]{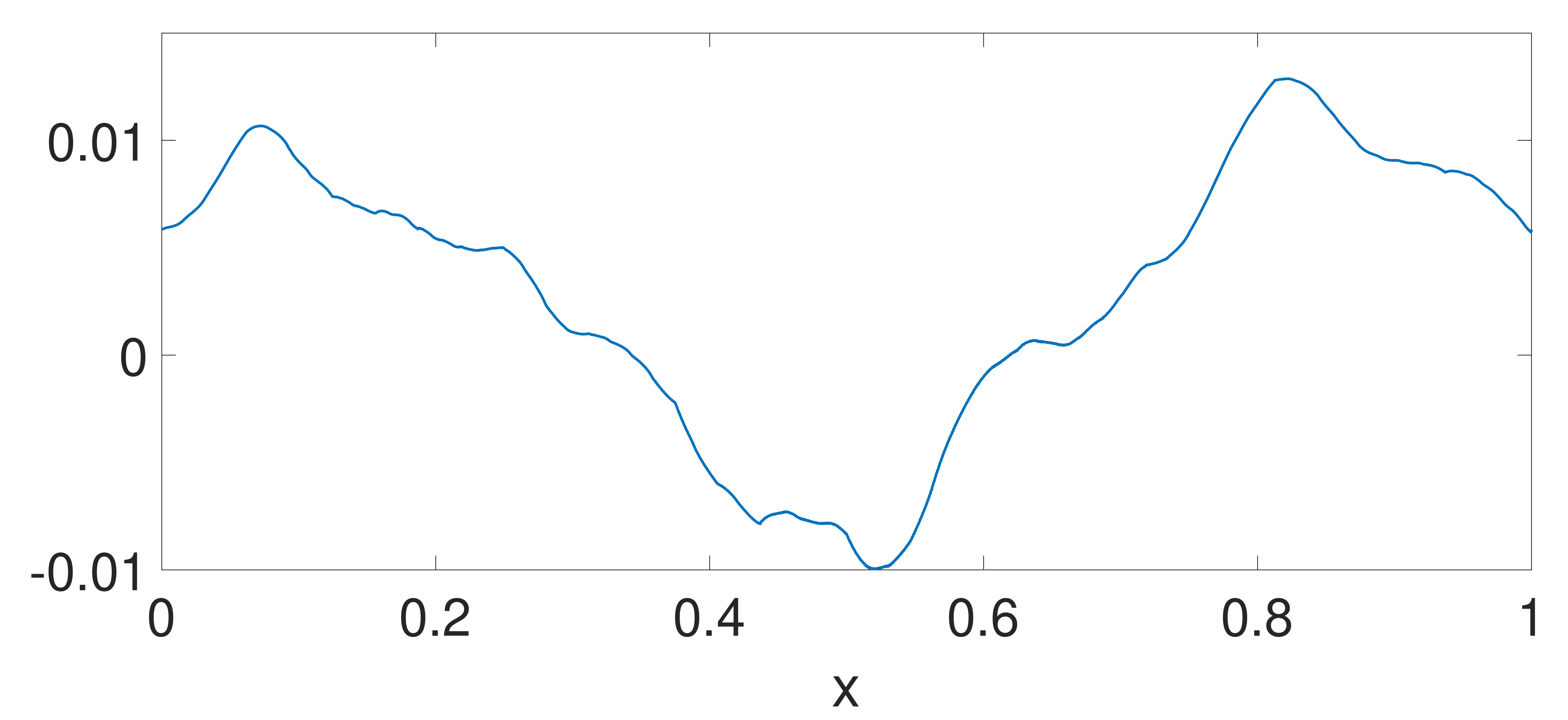}
        \caption{}
    \end{subfigure}    
    \caption{(a) An example of external potential $V$, predicted
    solution $u_{NN}$ and the corresponding reference solution
    $u_{\text{ref}}$ with $K=5$ and $\Ntrain=5000$. (b) Error between
    reference solution and predicted solution $u_{\text{ref}}-u_{NN}$.}
    \label{fig:LSE_train5000_K5}
\end{figure}

Table~\ref{tab:lse_K} presents results for different number of ReLU
layers $K$ with $\Ntrain = \Ntest = 5000$. As there are more ReLU layers,
we observe that the number of parameters increase monotonically and
both the train and test relative errors decrease monotonically, which
leads to a natural trade-off between accuracy and efficiency. According
to Table~\ref{tab:lse_K}, the performance improvement is marginal
beyond 5 ReLU layers. Figure~\ref{fig:LSE_train5000_K5}~(a) shows an
example of the external potential, the predicted and the corresponding
reference solution with $K=5$ and $\Ntrain=5000$. The first constraint
in \eqref{eq:linearschrodingereq} requires the norm of the discrete
solution $\vec{u}$ to be $\sqrt{2048}$ in our discretization settings. As
a result, we rescale the predicted solution $\calN_\theta$ to meet
the constraint. We notice that the NN result aligns well with the
reference solution, which implies that the solution map of the
linear eigenvalue problem can be trained under \stf{}. The error
between the reference solution and predicted solution is presented in
Figure~\ref{fig:LSE_train5000_K5}~(b).

\begin{table}[htp]
    \centering
    \begin{tabular}{ccccc}
        \toprule
        & Train rel err & Test rel err & Train energy & Test energy \\
        \toprule
        \ft{} & 3.66e-3 & 4.86e-3 & -112.51 & -110.77 \\
        \st{} & 7.11e-3 & 8.16e-3 & -113.31 & -111.60 \\
        \bottomrule
    \end{tabular}
    \caption{Relative error and predicted ground state energy $E$
    of $\calN_\theta$ trained under \stf{} and \ftf{} for linear
    Schr\"odinger equation with $K=5$ and $\Ntrain = \Ntest = 5000$.}
    \label{tab:lse_SNet}
\end{table}

We also compare the performance of different train frameworks,
as shown in Table~\ref{tab:lse_SNet}. The relative error of \ftf{}
is a little lower than that of \stf{} with the same number of ReLU
layers $K=5$ and train samples $\Ntrain=5000$. Such a difference in
relative error is mainly due to the different target of loss function.
Under \ftf{}, the loss function is least square between smallest
eigenvector and the NN output, which is consist with the relative
error defined as \eqref{eq:lschrodingerrelerr}. However, the loss
function under \stf{}, \eqref{eq:lschrodingerloss}, aims to minimize
the energy, which is inconsistent with the relative error. Hence
the difference between the relative errors for different train
frameworks is reasonable.  Predicted ground state energy of \stf{}
is lower than that of \ftf{}, which is also due to the different
loss function designs.  Considering the expensive data preparation
cost under the \ftf{}, i.e., solving \eqref{eq:linearschrodingereq}
for every input external potential $\{\vec{V}_i\}_{i=1}^{\Ntrain}$,
training under \stf{} is still desirable.

\subsection{Nonlinear Schr\"odinger equation}
\label{sec:nlschrodinger}

This section focuses on training the solution map of the smallest
eigenvalue problem of the one-dimensional nonlinear Schr\"odinger
equation (NLSE) as,
\begin{equation} \label{eq:nlse}
    \begin{split}
        & -\Delta u(x) + V(x) u(x) + \beta u(x)^3 = E u(x), \quad x
        \in \Omega = [0,1)\\
        & \mathrm{s.t.} \int_{\Omega} u(x)^2 \diff x = 1,~\text{and}
        \int_{\Omega} u(x) \diff x > 0,
    \end{split} 
\end{equation}
with periodic boundary condition.  The second positivity constraint in
\eqref{eq:nlse} can be dropped as in \eqref{eq:linearschrodingereq}
since the nonlinear term here is cubic. While, comparing to
\eqref{eq:linearschrodingereq}, the first constraint in \eqref{eq:nlse}
should be handled differently due to the nonlinearity and will be
taken care of in the NN design.  This NLSE \eqref{eq:nlse} is also
known as Gross-Pitaevskii (GP) equation in describing the single
particle properties of Bose-Eistein condensates. There is an associated
Gross-Pitaevskii energy functional,
\begin{equation} \label{eq:gpenergy}
    \calE[u(x)] = \innerproduct{\nabla u(x)}{\nabla
    u(x)} + \expectationvalue{V(x)}{u(x)} +
    \frac{\beta}{2}\innerproduct{u(x)}{u(x)^3},
\end{equation}
for positive $V(x)$ and $\beta$.  According to Theorem 2.1 in
\cite{Lieb1999}, the minimizer of the GP energy functional
\eqref{eq:gpenergy} is the eigenfunction of \eqref{eq:nlse}
corresponding to the smallest eigenvalue. 

In this section, the external potential is generated exactly the same
as that in \eqref{eq:potentialfun}, and then shifted such that the
minimum value of $V(x)$ equals to 1 in the observation of positivity
assumption on $V(x)$.  And $\beta$ here is set to be 10 such that
the problem is in the nonlinear regime. The NLSE \eqref{eq:nlse}
is discretized on a uniform grid of $2048$ points in the same way
as the linear Schr\"odinger equation \eqref{eq:linearschrodingereq}
in Section~\ref{sec:lschrodinger}.

While, the design of loss function for NLSE is more tricky.
Thanks to the GP energy functional, we define our loss function as
the discretized version of \eqref{eq:gpenergy},
\begin{equation} \label{eq:nlseloss}
    \ell \Big(\{\vec{V}_i\}_{i=1}^{\Ntrain}, \calA, \calN_\theta \Big)
    = \sum_{i=1}^{\Ntrain} \innerproduct{ \calN_\theta (\vec{V}_i)}{
    -\Delta \calN_\theta (\vec{V}_i)+\vec{V}_i\calN_\theta (\vec{V}_i)
    + \frac{\beta}{2} \calN_\theta (\vec{V}_i)^3},
\end{equation}
which again depends only on $\{\vec{V}_i\}_{i=1}^{\Ntrain}$, $\calA$,
$\calN_\theta$.

Through the derivative of \eqref{eq:nlseloss} with respect
to $\calN_\theta$ and train rule, minimizing our loss function
\eqref{eq:nlseloss} with respect to $\theta$ results the smallest
eigenvector $\calN_\theta(\vec{V}_i)$ for each $\vec{V}_i$ if the
NN $\calN_\theta$ is able to capture the solution map. However,
\eqref{eq:nlseloss} does not provide the smallest eigenvalue directly.
Instead, we calculate the smallest eigenvalue, i.e., the ground state
energy $E$, through a Rayleigh-quotient-like form as follows,
\begin{equation} \label{eq:nlseeigenvalue}
    E = \frac{\innerproduct{\calN_\theta (\vec{V})}{-\Delta
    \calN_\theta (\vec{V})+\vec{V}\calN_\theta (\vec{V}) +
    \beta \calN_\theta (\vec{V})^3}}{\innerproduct{\calN_\theta
    (\vec{V})}{\calN_\theta (\vec{V})}}.
\end{equation}

Similar as Section~\ref{sec:lschrodinger}, the loss function cannot
be used as a measure of the approximation accuracy of $\calN_\theta$.
We calculate the relative error on another set of $\Ntest$ random
external potential vectors of the same distribution as the train data,
$\{\vec{W}_i\}_{i=1}^{\Ntest}$. And the relative error is of the form,
\begin{equation} \label{eq:nlserelerr}
    \frac{1}{\Ntest} \sum_{i=1}^{\Ntest} \frac{\norm{ \vec{u}_i -
    \calN_\theta(\vec{W}_i)}}{\norm{ \vec{u}_i}},
\end{equation}
where $\vec{u}_i$ is the smallest eigenvector corresponding
to $\vec{W}_i$ for $i=1, \dots, \Ntest$.

We adopt the same $\calH$-net as in \cite{Fan2018} except that the
extra normalization layer here is
\begin{equation}
    \calN_\theta(\vec{V}) = \sqrt{2048}
    \frac{\widetilde{\calN}_\theta}{\bigl\lVert\widetilde{\calN}_\theta \bigr\rVert},
\end{equation}
where $\widetilde{\calN}_\theta$ is the regular
$\calH$-net~\cite{Fan2018} and $\calN_\theta$ is the NN used in this
section. This extra layer makes sure the norm of the NN output equals
to $\sqrt{2048}$, which agrees with the discretized version of the
first constraint in \eqref{eq:nlse} and also agrees with the reference
solution generated through the traditional method~\cite{Bao2004},
which is under the same settings as \cite{Fan2018}. All training
details are the same as that in Section~\ref{sec:lschrodinger}.

\subsubsection*{Numerical results}

We first compare the performance of $\calN_\theta$ trained under \stf{}
for different number of train data set size $\Ntrain$ and different
number of ReLU layers $K$ through numerical experiments.

\begin{table}[htp]
    \centering
    \begin{tabular}{cccc}
        \toprule
        $\Ntrain$  & $\Ntest$ & Train rel err  & Test rel err \\
        \toprule
        500        & 5000     & 2.52e-2        & 2.69e-2 \\
        1000       & 5000     & 3.01e-2        & 3.18e-2 \\
        5000       & 5000     & 4.98e-3        & 5.28e-3 \\
        20000      & 20000    & 5.24e-3        & 5.25e-3  \\
        \bottomrule
    \end{tabular}
    \caption{Relative error of $\calN_\theta$ with $K=5$ trained under
    \stf{} for NLSE given different sizes of train and test data set.}
    \label{tab:nlse_size}
\end{table}

\begin{table}[htp]
    \centering
    \begin{tabular}{cccc}
        \toprule
        $K$ & $N_{params}$ & Train rel err & Test rel err \\
        \toprule
        1   & 15184        & 1.64e-1       & 1.65e-1 \\
        3   & 34236        & 1.49e-2       & 1.51e-2 \\
        5   & 57156        & 4.98e-3       & 5.28e-3 \\
        7   & 83944        & 3.65e-3       & 3.95e-3 \\
        \bottomrule
    \end{tabular}
    \caption{Relative error of $\calN_\theta$ trained under \stf{}
    for NLSE with different number of ReLU layers $K$. The train and
    test data sets are of size $\Ntrain=5000$ and $\Ntest=5000$.}
    \label{tab:nlse_K}
\end{table}

\begin{figure}[htp]\centering
    \begin{subfigure}{0.495\textwidth}
        \includegraphics[width=\hsize]{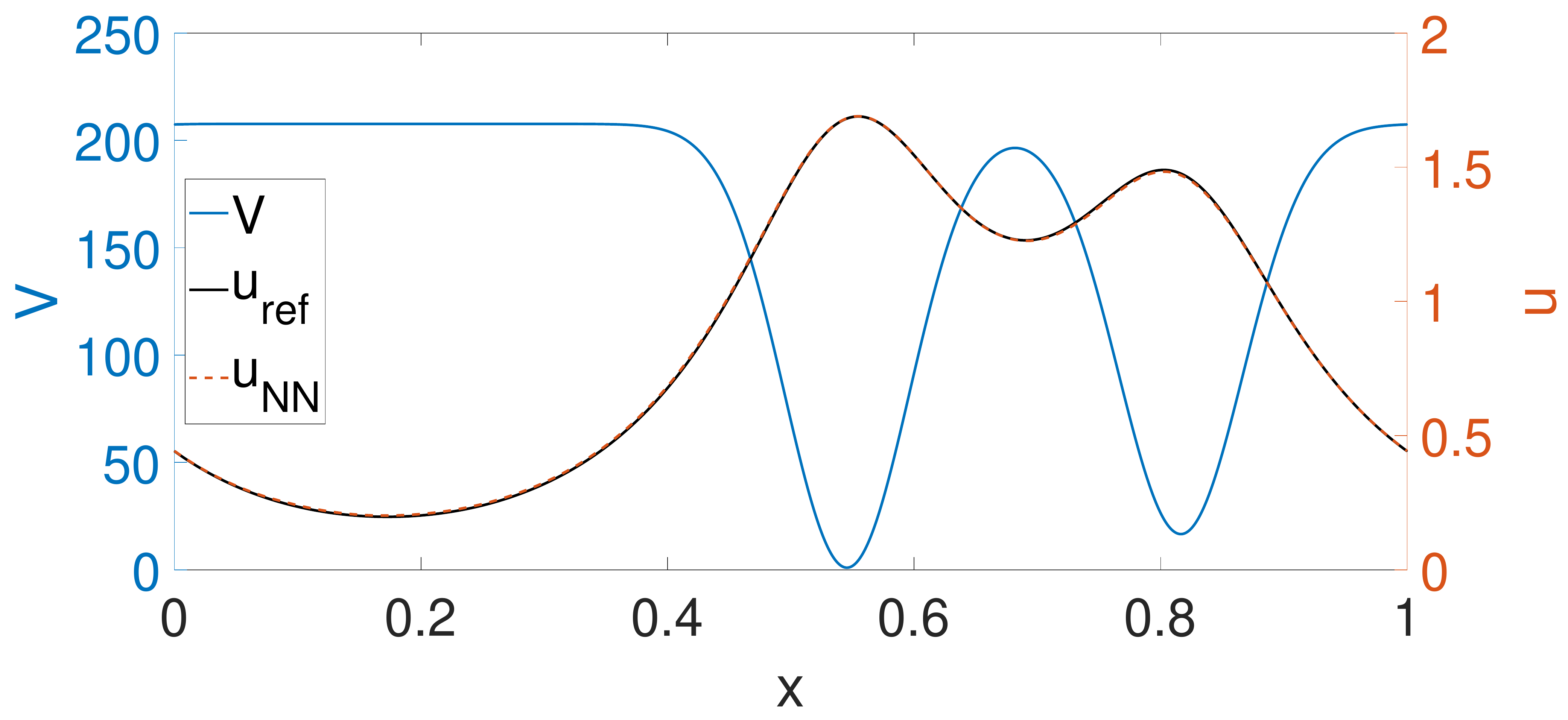}
        \caption{}
    \end{subfigure}    
    \begin{subfigure}{0.495\textwidth}
        \includegraphics[width=\hsize]{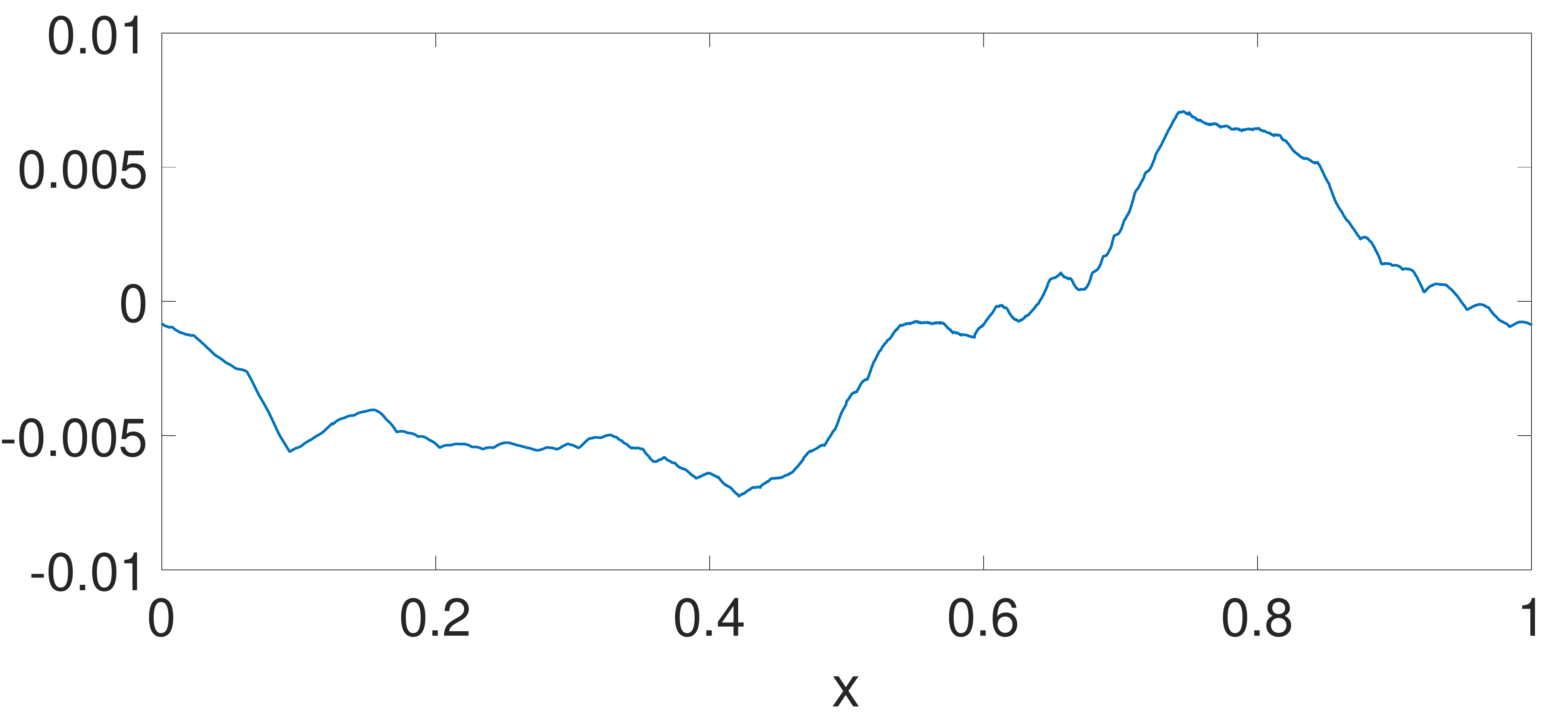}
        \caption{}
    \end{subfigure}    
    \caption{(a) An example of external potential $V$, predicted
    solution $u_{NN}$ and the corresponding reference solution
    $u_{\text{ref}}$ with $K=5$ and $\Ntrain=5000$. (b)
    Error between reference solution and predicted solution
    $u_{\text{ref}}-u_{NN}$.} \label{fig:NLSE_train5000_K5}
\end{figure}

Table~\ref{tab:nlse_size} and Table~\ref{tab:nlse_K} present results
for different $\Ntrain$, $\Ntest$, and for different number of ReLU
layers $K$ respectively. Similar as in Section~\ref{sec:lschrodinger},
the test relative error decreases as $\Ntrain$ and $K$
increases. The results are near optimal with $\Ntrain=5000$ and
$K=5$. Figure~\ref{fig:NLSE_train5000_K5}~(a) shows an example of the
external potential, the predicted and the corresponding reference
solution. The error between the reference solution and predicted
solution is presented in Figure~\ref{fig:NLSE_train5000_K5}~(b).
All comments in Section~\ref{sec:lschrodinger} apply here.

\begin{table}[htp]
    \centering
    \begin{tabular}{ccccc}
        \toprule
        & Train rel err & Test rel err & Train energy & Test energy \\
        \toprule
        \ft{} & 1.68e-3 & 2.02e-3 & 152.43 & 152.94 \\
        \st{} & 4.98e-3 & 5.28e-3 & 152.34 & 152.84 \\
        \bottomrule
    \end{tabular}
    \caption{Relative error and predicted ground state energy $E$
    of $\calN_\theta$ trained under \stf{} and \ftf{} for NLSE with
    $K=5$ and $\Ntrain = \Ntest = 5000$.} \label{tab:nlse_SNet}
\end{table}

We also compare the performance of different train frameworks,
as shown in Table~\ref{tab:nlse_SNet}. Again similar as in
Section~\ref{sec:lschrodinger}, the relative error of \ftf{} is
a little lower than that of \stf{} with the same number of ReLU
layers $K=5$ and train samples $\Ntrain=5000$. However, predicted
ground state energy of \stf{} is lower than that of \ftf{} due to the
different choices of loss functions. We observe nonlinear behavior
for the solution near $x=0$ and the approximation error is even
smaller than that in Section~\ref{sec:lschrodinger}. Hence training
the $\calN_\theta$ proposed in \cite{Fan2018} under \stf{} is able
to achieve similar accuracy with similar training computational cost
but saves the expensive train data preparation step comparing against
training under \ftf{}.

\section{Conclusion}
\label{sec:conclusion}

We propose a novel training framework named \stf{} to train NN
in representing low dimensional solution maps of physical models.
Since physical models have fixed forward maps usually in the form
of PDEs, NN can be viewed as the ansatz of the solution map and be
trained variationally with unlabeled input functions through a loss
function containing forward maps, i.e., $\ell(\{\vec{f}_i\}, \calA,
\calN_\theta)$ for $\{\vec{f}_i\}$, $\calA$, and $\calN_\theta$
being the input functions, forward map, and NN respectively. Training
under \stf{} is able to avoid the expensive data preparation step,
which prepares labels for input functions through costly traditional
solvers, and still captures the solution map adapted to the input
data distribution.

The power of \stf{} is illustrated through four examples, solving
linear and nonlinear elliptic equations and solving the ground
state of linear and nonlinear Schr\"odinger equations. For linear
elliptic equations, we use $\calH$-matrix structure as the ansatz
and train via the loss as \eqref{eq:loss}. The trained solution
map outperforms the traditional $\calH$-matrix obtained from SVD
truncation. For nonlinear elliptic equations, we use one wide fully
connected layer using ReLU activation NN as the ansatz and train
via the same loss. Without labeling the input data, \stf{} is able to
achieve the solution map adaptive to the input distribution whereas the
traditional training framework fails in training. Finally, for both
linear and nonlinear Schr\"odinger equations, we adopt variational
representation of the ground state energy as the loss function and
train $\calH$-nets~\cite{Fan2018} under \stf{}. Lower ground state
energy is obtained via \stf{} comparing to the traditional \ftf{}.

\medskip
\noindent
{\bf Acknowledgments.}  This work is partially
supported by the National Science Foundation under awards OAC-1450280
and DMS-1454939.

\bibliographystyle{apalike}
\bibliography{library}

\begin{appendices}

    \section{$\calH$-matrix structure}
    \label{sec:append_hmatrix}
    
    Assume the mapping between input vector $f$ and output vector $u$
    is a matrix $A\in\mathbb{R}^{N^2 \times N^2}$, where $N=2^L m$
    and admits the two dimensional $\calH$-matrix structure. In order
    to simplify the description below, we introduce a few handy
    notations. The bracket of an integer is adopted to denote the
    set of nonnegative integers smaller than the given one, i.e., $[n] =
    \{0, 1, \dots, n-1\}$. Although $u$ and $f$ are always viewed as
    vectors, they are functions on a two dimensional grid. To avoid
    complicated notations, we adopt two input indices for them.
    Further a fully connected layer (dense layer) with input size
    $n_1$ and output size $n_2$ is denoted as $\calD_{n_1}^{n_2}$ and
    the ReLU activation function is denoted as $\sigma(\cdot)$. The
    corresponding $\calH$-matrix neural network structure with rank
    $r$ can be constructed as follows.
    \begin{itemize}
        \item {\bf Level $\ell = 1, 2, \dots, L$.} On level $\ell$,
            the indices are split into $2^{\ell}$ parts, denoted as
            $\calI^{\ell}_{i} = 2^{L-\ell}m \cdot i + [2^{L-\ell}m]$
            for $i \in [2^\ell]$.  Hence vector $u \big( \calI^\ell_a,
            \calI^\ell_b \big)$ and $f \big( \calI^\ell_i, \calI^\ell_j
            \big)$ are of length $k = 2^{2L-2\ell}m^2$ for any $a,b,i,j
            \in [2^\ell]$.  Then the operation on level $\ell$ is
            defined as,
            \begin{equation} \label{eq:hmatlr}
                u \big( \calI^{\ell}_a, \calI^{\ell}_b \big) =
                u \big( \calI^{\ell}_a, \calI^{\ell}_b \big) +
                \sum_{\substack{i,j \in [2^\ell] \\ i \neq a,
                j \neq b}} \calD^k_r \sigma \Big( \calD^r_k
                f\big(\calI^{\ell}_i, \calI^{\ell}_j \big) \Big)
            \end{equation}
            for $a, b \in [2^\ell]$. When the activation function
            $\sigma$ is removed from \eqref{eq:hmatlr}, the operation in
            the summation is a low-rank factorization of the mapping
            between grid $\big( \calI^\ell_i, \calI^\ell_j \big)$ and
            $\big( \calI^\ell_a, \calI^\ell_b \big)$, which is known as
            the far-field interaction in $\calH$-matrix literature.

        \item {\bf Diagonal Level.} On this level, the same indices as
            on level $L$ are used, i.e., $\calI^{L}_{i} = m \cdot i
            + [m]$ for $i \in [2^L]$. The operation on this level is
            defined as,
            \begin{equation} \label{eq:hmatdense}
                u \big( \calI^L_a, \calI^L_b \big) =
                u \big( \calI^L_a, \calI^L_b \big) +
                \calD^{m^2}_{m^2}
                f \big( \calI^L_a, \calI^L_b \big)
            \end{equation}
            for $a, b \in [2^L]$. This operation is known as the
            near-field (local) interaction in $\calH$-matrix
            literature.
    \end{itemize}
    
    The description of $\calH$-matrix is abstract and lacks the
    domain decomposition intuition behind it. Readers are referred
    to \cite{Hackbusch1999} for more details about $\calH$-matrix.

    In Section~\ref{sec:solvelinearsystem}, NN-$\calH$-matrix refers
    to the neural network without activation function $\sigma$ whereas
    NLNN-$\calH$-matrix refers to the neural network with $\sigma$.
    Further, when SVD initialization is used, we first calculate
    the rank $r$ truncated SVD of the submatrix of $A$ mapping from
    $\big( \calI^\ell_i, \calI^\ell_j \big)$ to $\big( \calI^\ell_a,
    \calI^\ell_b \big)$ and then initialize $\calD^k_r$ by the product
    of left singular vectors and singular values, and $\calD^r_k$
    by the right singular vectors.
    
    \section{$\calH$-net~\cite{Fan2018} structure}
    \label{sec:append_hnet}
    
    Since $\calH$-net is applied to one dimensional problems, we
    introduce the one dimensional version here. Assume that both the
    input vector $f$ and the output vector $u$ are discretized on $N$
    points, where $N=2^L m$. We follow the notations in~\cite{Fan2018}
    and assume that all the relevant tensors will be appropriately
    reshaped or padded for simplicity. A tensor $\xi$ of size
    $\alpha \times N_x$ is connected to a tensor $\zeta$ of size
    $\alpha'\times N_x'$ by a locally connected (LC) network if
    \begin{equation}
       \zeta_{c',i} = \phi \left( \sum_{j=(i-1)s+1}^{(i-1)s+w}
       \sum_{c=1}^{\alpha} W_{c',c;i,j} \xi_{c,j} + b_{c',i}\right),
       \quad i=1,\ldots,N_{x}',~c'=1,\ldots,\alpha',
    \end{equation}
    where $w$ is the kernel window size and $s$ is the
    stride. Three kinds of LC networks are combined in
    $\calH$-net. Among them, $\text{LCR}[\phi;N_x,N_x',\alpha']$
    denotes the restriction network where $s=w=\frac{N_x}{N_x'}$
    and $\alpha=1$, $\text{LCK}[\phi;N_x,\alpha,\alpha',w]$
    denotes the kernel network where $s=1$ and $N_x'=N_x$, and
    $\text{LCI}[\phi;N_x,\alpha,\alpha']$ denotes the interpolation
    network where $s=w=1$ and $N_x'=N_x$. The ReLU activation function
    is denoted as $\sigma(\cdot)$. The corresponding $\calH$-net
    structure with rank $r$ can be constructed as follows.
    \begin{itemize}
    	\item {\bf Level $\ell = 2, 3, \dots, L$.} On level $\ell$,
	    the indices are split into $2^{\ell}$ parts. The operation
	    on level $\ell$ is defined as,
    	\begin{equation}
    	    \begin{split}
    	    & \xi_0 = \text{LCR}[\text{linear};N,2^{\ell},r](f)\\
    	    & \xi_k = \text{LCK}[\sigma;2^{\ell},r,r,2n_b^{({\ell})}+1](\xi_{k-1}), \quad k=1,\ldots,K\\
    	    & u = u+\text{LCI}[\text{linear};2^{\ell},r,\frac{N}{2^{\ell}}](\xi_K)
    	    \end{split} 
    	\end{equation}
    	where $n_b^{({\ell})}$ is 2 for $\ell=2$ and 3 for $\ell\geq 3$.
    
	\item {\bf Diagonal Level.} On this level, the same indices
	as on level $L$ are used. The operation on this level is
	defined as,
    	\begin{equation}
    	\begin{split}
    	& \xi_k = \text{LCK}[\sigma;2^L,m,m,2n_b^{(\text{ad})}+1](\xi_{k-1}), \quad k=1,\ldots,K-1\\
    	& u = u+\text{LCK}[\text{linear};2^L,m,m,2n_b^{(\text{ad})}+1](\xi_{K-1})
    	\end{split} 
    	\end{equation}
    	where $\xi_0 = f$ and $n_b^{(\text{ad})}=1$.
    \end{itemize}
    
    Readers are referred to~\cite{Fan2018} for the intuition behind
    the design of $\calH$-net structure.

\end{appendices}

\end{document}